\documentclass{amsart}

\usepackage{a4,amsfonts,amssymb,amsmath,amscd,xy}
\title{Double Poisson Structures on Finite Dimensional Semi-Simple Algebras}
\author{Geert Van de Weyer
}
\address{Department of Mathematics and Computer Science\\ University of Antwerp\\ B-2020 Antwerpen (Belgium)}
\email{geert.vandeweyer@ua.ac.be}
\thanks{The author is postdoctoral fellow of the Fund for Scientific Research \--- Flanders (F.W.O.-Vlaanderen)(Belgium).}

\newcommand{\CC}{\mathbb{C}}

\newcommand{\ZZ}{\mathbb{Z}}

\newcommand{\Der}{\mathrm{Der}}
\newcommand{\DDer}{\mathbb{D}\mathrm{er}}
\newcommand{\vtx}[1]{*+[o][F-]{\scriptscriptstyle #1}}

\newcommand{\DDQ}{\mathbb{D}S}
\newcommand{\rep}{\mathsf{rep}}
\newcommand{\iss}{\mathsf{iss}}
\newcommand{\GL}{\mathsf{GL}}
\newcommand{\SL}{\mathsf{SL}}
\newcommand{\db}[1]{\{\!\!\{ #1 \}\!\!\}}

\newtheorem{theorem}{Theorem}
\newtheorem{prop}{Proposition}
\newtheorem{lemma}{Lemma}
\newtheorem{cor}{Corollary}
\newtheorem{defn}{Definition}

\xyoption{all}
\allowdisplaybreaks
\begin{document}

\bibliographystyle{plain}
\maketitle

\begin{abstract}
We give a description of the bimodule of double derivations $\DDer(S)$ of a finite dimensional semi-simple algebra $S$ and its double Schouten bracket in terms of a quiver. This description is used to determine which degree two monomials in $T_S\DDer(S)$ induce double Poisson brackets on $S$. In case $S = \CC^{\oplus n}$, a criterion for any degree two element to give a double Poisson bracket is deduced. For $S = \CC^{\oplus n}$ and $S' = \CC^{\oplus m}$ the induced Poisson bracket on the variety of isomorphism classes of semi-simple representations $\iss_n(S*T)$ of the free product $S*T$ is given.
\end{abstract}

\section{Introduction}
Throughout this note, we will work over an algebraically closed field of characteristic $0$ which we denote by $\CC$. Unadorned tensor products will be over the base field $\CC$. 

Double Poisson algebras were introduced by M.~Van~den~Bergh in \cite{MichelDPA} as a generalization of classical Poisson geometry to the setting of noncommutative geometry. The key fact being that an algebra $A$ equipped with a double Poisson bracket has a canonical Poisson structure on all its finite dimensional representation spaces $\rep_n(A)$. More specifically, a double Poisson algebra $A$ is an associative unital algebra equipped with a linear map
$$\db{-,-}:A\otimes A\rightarrow A\otimes A$$
that is a derivation in its second argument for the outer $A$-bimodule structure on $A\otimes A$, where the outer action of $A$ on $A\otimes A$ is defined as $a.a'\otimes a''.b := (aa')\otimes (a''b)$. Furthermore, we must have that $\db{a,b} = -\db{b,a}^o$ and that the double Jacobi identity holds for all $a,b,c\in A$:
\begin{align*}
&\db{a,\db{b,c}'}\otimes\db{b,c}'' + \db{c,a}''\otimes \db{b,\db{c,a}'}\\ & + \db{c,\db{a,b}'}''\otimes\db{a,b}''\otimes\db{c,\db{a,b}'}'
= 0,
\end{align*}
where we used Sweedler notation, that is $\db{x,y} = \sum \db{x,y}'\otimes\db{x,y}''$ for all $x,y\in A$. Such a map is called a \emph{double Poisson bracket}. 

A double Poisson bracket yields, for each $n$, a classical Poisson bracket on the coordinate ring $\CC[\rep_n(A)]$ of the variety of $n$-dimensional representations of $A$ through $\{a_{ij},b_{k\ell}\} := \db{a,b}_{kj}'\db{a,b}_{i\ell}$. This bracket restricts to a Poisson bracket on $\CC[\rep_n(A)]^{\GL_n}$, the coordinate ring of the quotient variety $\iss_n(A)$ under the action of the natural symmetry group $\GL_n$ of $\rep_n(A)$.

In this paper, we study double Poisson brackets on a direct sum $S = M_{d_1}(\CC)\oplus \dots \oplus M_{d_k}(\CC)$ of matrix algebras over $\CC$. Because such algebras are smooth, we know from \cite{MichelDPA} that all double Poisson brackets are determined by double Poisson tensors. That is, elements of degree $2$ in $\DDQ = T_S\DDer(S)$ where $\DDer(S) = \Der(S,S\otimes S)$ is the module of double derivations. That is, the module of derivations from $S$ to the $S$-bimodule $S\otimes S$, where the $S$-action on $S\otimes S$ is the outer action. $\DDer(S)$ is an $S$-bimodule through the inner action: $(s.\vartheta.t)(u) = \vartheta(u)'t\otimes s\vartheta(u)''$. A first important result is the explicit description of $\DDer(S)$ and $\DDer_T(S)$. Here, $\DDer_T(S)$ is the bimodule of $T$-linear double derivations with $T\subset S$ a subalgebra. That is, double derivations that are identically zero on $T$. We have that
$$\DDer(S) \cong \bigoplus_{i=1}^k M_{d_i}(\CC)^{\oplus d_i^2-1} \oplus
\bigoplus_{i\neq j} M_{d_i\times d_j}(\CC)^{\oplus d_id_j}$$
as $S$-bimodules where $S$ acts on the right hand side expression by matrix multiplication. If
$T=M_{e_1}(\CC)\oplus\dots M_{e_\ell}(\CC)$ is a finite dimensional semi-simple subalgebra of $S$ with Bratelli diagram with respect to $S$ given by $(a_{ij})_{(i,j)=(1,1)}^{(k,\ell)}$, then
$$\DDer_T(S) \cong \bigoplus_{i=1}^k M_{d_i}(\CC)^{\oplus r_i} \oplus
\bigoplus_{i\neq j} M_{d_i\times d_j}(\CC)^{\oplus r_{ij}}$$
as $S$-bimodules, with $r_i = \sum_{u=1}^l a_{iu}^2 - 1$ and $r_{ij} = \sum_{u=1}^l a_{iu}a_{ju}$.
These descriptions are formulated in Section \ref{DDQ}, Theorem \ref{DoubleDerivationsNotRelative} and Theorem \ref{DoubleDerivationsRelative}.

Using these two theorems, we are able to formulate, in Theorem \ref{DDQGradedLie}, an explicit description of the graded Lie algebra $\DDQ/[\DDQ,\DDQ][1]$, where the bracket on $\DDQ/[\DDQ,\DDQ][1]$ is the bracket associated to the double Schouten-Nijenhuis bracket on $\DDQ$. This description is formulated in terms of the \emph{double derivation quiver} $Q_S$ associated to $S$. Assign to $S$ a quiver $\overline{Q}_S$ on $k$ vertices with $d_id_j$ arrows between each two vertices $i\neq j$ and $d_i^2-1$ loops in all vertices $i$, where the arrows are indexed by index sets $C_{ji} = \{1,\dots,d_j\}\times\{1,\dots, d_i\}$ if $i\neq j$ and $C_{ii} = \{1,\dots,d_j\}\times\{1,\dots, d_i\}\backslash\{(1,1)\}$. Then
$\DDQ/[\DDQ,\DDQ][1]$ is isomorphic as a graded Lie algebra to $\CC Q_S/[\CC Q_S,\CC Q_S]_{super}$, where the bracket on two words $\omega_1 = v_1\dots v_n$ and $\omega_2 = u_1\dots u_m$ in $\CC\overline{Q}_S/[\CC Q_S,\CC Q_S]_{super}$ is depicted in Figure \ref{GradedBracket}.

\begin{figure}
$$
\xymatrix@R=1.75pc@C=1.75pc{
\save[0,1].[7,1]!C*+<10pt,0pt>\frm{(}\restore
\save[0,1].[7,1]!C+<-55pt,0pt>
\drop{\scriptstyle\sum\limits_{a\in (\overline{Q}_S)_0} (-1)^{(i+j)(n-1)}}
\restore
\save[0,4].[7,4]*+<10pt,0pt>\frm{)}\restore
\save[0,7].[7,7]*+<10pt,0pt>\frm{(}\restore
\save[0,7].[7,7]!C+<-42pt,0pt>
\drop{\scriptstyle-(-1)^{(i+j+1)(n-1)}}
\restore
\save[0,10].[7,10]*+<10pt,0pt>\frm{)}\restore
&
&\bullet \ar[r]^{u_1} & \bullet \ar[dr]^{u_2} 
&&& 
& &\bullet \ar[r]^{u_1} & \bullet \ar[dr]^{u_2} & 
\\ & 
\bullet\ar[ur]^{u_m}&& & \bullet \ar[d]^{\vdots}
&&&
\bullet\ar[ur]^{u_m}&& & \bullet \ar[d]^{\vdots}
\\ &
\bullet\ar[u]\ar@{}[urrr]|{\omega_2}& & & \bullet \ar[dl]^{u_i}
&&&
\bullet\ar[u]\ar@{}[urrr]|{\omega_2}&& & \bullet \ar[dl]^{u_i}
\\ &
 &\bullet\ar[ul]^{u_{i+2}}& a \ar[l]^{(p,q)}|{~~//}\ar@/^.75pc/@{.}[d]|=&
&&&
&a \ar[ul]^{u_{i+2}} & \bullet\ar[l]^{(p,q)}|{~~//}\ar[d]^{(r,q)} &
 \\ &
 &\bullet \ar[u]^{(p,s)}\ar[r]^{(q,s)}|{//~~}  & a\ar[dr]^{v_{j+2}}& &
&& 
& a\ar[r]^{(r,p)}|{//~~}\ar@/^.5pc/@{.}[u]|{=} & \bullet\ar[dr]^{v_{j+2}}& 
\\   &
\bullet\ar[ur]^{v_j}&& & \bullet \ar[d]^{\vdots}
&& & 
\bullet\ar[ur]^{v_j}&& & \bullet \ar[d]^{\vdots}
\\ &
\bullet\ar[u]^{\vdots}\ar@{}[urrr]|{\omega_1}&& & \bullet \ar[dl]^{v_n}
&&& 
\bullet\ar[u]^{\vdots}\ar@{}[urrr]|{\omega_1}&& & \bullet \ar[dl]^{v_n}
\\ &
&\bullet\ar[ul]^{v_{2}}& \bullet \ar[l]^{v_1}
&&& 
& &\bullet\ar[ul]^{v_{2}}& \bullet \ar[l]^{v_1} & & 
\\
}
$$
\caption{\label{GradedBracket} The graded Lie bracket on $\CC Q/[\CC Q,\CC Q]_{super}$.}
\end{figure}
This result can then be used to dermine all monomials  of degree $2$ in $\DDQ$ that yield nontrivial double Poisson structures on $S$. In Section \ref{DPLC} we use this result to compute the first double Poisson-Lychnerowicz cohomology groups for $S$.

Although the representation varieties of finite dimensional semi-simple algebras are rather simple and the quotient varieties consist of a finite number of points, double Poisson structures on these algebras yield interesting noncommutative geometry as they can be extended to double Poisson structures on the free product of such algebras. For such a free product $S*T$, the quotient variety $\iss_n(S*T)$ is no longer trivial and double Poisson structures can yield nontrivial Poisson structures on this variety. In the final section of this paper, we give an explicit description of the Poisson brackets on the quotient variety $\iss_n(\CC^{\oplus p}*\CC^{\oplus q})$.

\noindent\textbf{Acknowledgements.} Part of the work presented here was done while the author was visiting the University of Leeds, supported by a travel grant from the Fund for Scientific Research \--- Flanders (F.W.O.-Vlaanderen)(Belgium). The author would like to thank Peter J{\o}rgensen for inviting him and Bill Crawley-Boevey for several interesting discussions.

\section{Double Derivations and Quivers}\label{DDQ}
Consider a finite dimensional semi simple algebra $S = M_{d_1}(\CC)\oplus \dots \oplus M_{d_k}(\CC)$. As seen in the introduction, the set of double derivations $\DDer(S) := \Der(S,S\otimes S)$ for the outer action of $S$ on $S\otimes S$ can be equipped with a $S$-bimodule structure using the inner action of $S$ on $S\otimes S$, defined as $a\circ a'\otimes a''\circ b := (a'b)\otimes (aa'')$. The $S$-bimodule of double derivations can be described as
\begin{theorem}\label{DoubleDerivationsNotRelative}
$$\DDer(S) \cong \bigoplus_{i=1}^k M_{d_i}(\CC)^{\oplus d_i^2-1} \oplus
\bigoplus_{i\neq j} M_{d_i\times d_j}(\CC)^{\oplus d_id_j}$$
as $S$-modules where the actions on $\Der(S,S\otimes S)$ are the inner actions and the actions on the right hand side are just matrix multiplication:
\begin{align*}
& (s_1,\dots,s_k).(m_1,\dots,m_k,m_{12},\dots, m_{ij},\dots, m_{(k-1)k},m_{k(k-1)}).(t_1,\dots,t_k) \\
& = (s_1m_1t_1,\dots,s_km_kt_k,s_1m_{12}t_2,\dots,s_im_{ij}t_j,\dots,s_{k-1}m_{(k-1)k}t_k,s_km_{k(k-1)}t_{k-1}).
\end{align*}
\end{theorem}
\begin{proof}
First of all note that because $S$ is semi simple all derivations must be inner, that is, any derivation $d:S\rightarrow S\otimes S$ is of the form
$$d(y) = d_x(y) = x'\otimes (x''y) - (yx')\otimes x''$$
for some $x = x'\otimes x'' \in S\otimes S$.
So we consider the morphism 
$$D:S\otimes S\twoheadrightarrow \DDer(S) : x\mapsto d_x$$
of $S\otimes S$ modules for the inner action and will compute its kernel.

First of all note that $S\otimes S$ for the inner action is isomorphic to $\bigoplus_{i=1}^k M_{d_i}(\CC)^{\oplus d_i^2} \oplus
\bigoplus_{i\neq j} M_{d_i\otimes d_j}(\CC)^{\oplus d_id_j}$. Indeed, denote by $e_{pq}^r$ the standard basis elements for $M_{d_r}(\CC)$, then for $s\in M_{d_j}(\CC)$ and $t\in M_{d_i}(\CC)$ we have
\begin{eqnarray*}
s.(e^i_{pq}\otimes e^j_{rs}).t &=& (e^i_{pq}t)\otimes (se^j_{rs}) \\
& = & \sum_{a=1}^{d_j}\sum_{b=1}^{d_i} e^i_{pb} t_{qb}\otimes e^j_{as} s_{ar}.
\end{eqnarray*}
That is, on the subspace spanned by
$$f^{ps}_{ab} := e^i_{pb}\otimes e^j_{as},~1\leq a\leq d_j,~1\leq b\leq d_i$$
we have
$$s.f^{ps}_{rq}.t = \sum_{a=1}^{d_j}\sum_{b=1}^{d_i}  s_{ar}f^{ps}_{ab}t_{qb}.$$
But this is nothing else than taking the $n\times m$ matrix with $1$ on row $r$ and column $q$ and zeroes elsewhere, and multiplying it on the left by the matrix $s$ and on the right by the matrix $t$, and we have $d_id_j$ such subspaces for each $i$ and $j$.

In order to compute the kernel of the map $D$, decompose $x\in S\otimes S$ as $x = (x^{p}, x^{qr})_{1\leq p \leq k, 1\leq q\neq r \leq k}$ with $x^p \in M_{d_p}(\CC)\otimes M_{d_p}(\CC)$ and $x^{qr}\in M_{d_q}(\CC)\otimes M_{d_r}(\CC)$.
For this decomposition we have
$$
d_x(y_1,\dots,y_k) = (x^py_p - y_px^p, x^{qr}y_r - y_qx^{qr})_{1\leq p \leq k, 1\leq q\neq r \leq k}.
$$
Now if $d_x = 0$, fix $r$ and take $y_r = 1$ and $y_s = 0$ for $s\neq r$ which yields $x^{qr} = 0$ for all $q$. Taking $y_r = e^r_{ii}$ yields $x^re^r_{ii} - e^r_{ii}x^r = 0$, whence
$$x^r = \sum_{a,b,p=1}^{d_r} x^r_{abp} f^{pp}_{ab}.$$
Finally, taking $y_r = e^r_{ij}$ yields $x^r e^r_{ij} - e^r_{ij} x^r = 0$, which gives $x_{abp} = x_{abq} = x_{ab}$ for all $1\leq p,q,a,b \leq d_r$. But then
$$x^r = \sum_{a,b,p=1}^{d_r} x^r_{ab} f^{pp}_{ab}.$$
So for the decomposition of $S\otimes S$ as a direct sum of modules of matrices, we get $Ker(D) = M_{d_1}(\CC)\oplus\dots\oplus M_{d_k}(\CC)$ from which the theorem follows.
\end{proof}
This description can be encoded using the language of quivers.
\begin{defn}\label{DDQuiver}
Let $S$ be a finite dimensional semi-simple algebra $M_{d_1}(\CC)\oplus \dots \oplus M_{d_k}(\CC)$, then the \emph{double derivation quiver $Q_S$} for $S$ is the quiver with $k$ vertices, $d_i^2-1$ loops in vertex $i$ and $d_id_j$ arrows from vertex $i$ to vertex $j$ for all $i\neq j$. The arrows $j\rightarrow i$, $i\neq j$ are indexed by a bi-index in the \emph{colour index set} $C_{ji} = \{1,\dots,d_j\}\times\{1,\dots,d_i\}$. The loops are indexed by a bi-index in the set $C_{ii} = \{1,\dots, d_i\}\times \{1,\dots, d_i\}\backslash\{(1,1)\}$. The first index will be called the \emph{primary colour} of the arrow and the second index will be called the \emph{secondary colour} of the arrow.
\end{defn}
As we already indicated in the introduction, the tensor algebra $\DDQ = T_S\Der(S,S\otimes S)$ plays an important role in the study of double Poisson brackets. This algebra can be seen as a 'path algebra' of the double derivation quiver $Q_S$ of $S$. In order to clarify what we mean, note first of all that $\DDer(S)$ is generated as an $S$-module by the elements
$$x^i_{pq} = e^i_{p1}\otimes e^i_{1q} \in M_{d_i}(\CC)$$
with $(p,q)\neq (1,1)$ and
$$y^{ij}_{pq} = e^i_{p1}\otimes e^j_{1q} \in M_{d_i}(\CC)$$
with $i\neq j$. Identifying $x^i_{pq}$ with the loops in vertex $i$ of $Q_S$ and $y^{ij}_{pq}$ with the arrows from vertex $j$ to vertex $i$, the multiplication of two elements of degree $1$ in $D_S$, say $u = sat$ and $v = s'a't'$ with $a$ and $a'$ arrows in $Q_S$ is then easily seen to be zero unless $h(a) = t(a')$.

We can easily see that in case $S = \CC^{\oplus n}$ the algebra $\DDQ$ is indeed equal to the path algebra of the double derivation quiver.
\begin{prop}
Let $S = \CC^{\oplus n}$, then $\DDQ = \CC Q_S$.
\end{prop}

Recall that for a subalgebra $T\subset S$, a double derivation $\vartheta$ is called $T$-linear if and only if $\vartheta(T) = 0$. Next, recall that a finite dimensional semi-simple subalgebra $T=M_{e_1}(\CC)\oplus\dots M_{e_\ell}(\CC)$ of $S$ is completely determined by its Bratelli diagram $(a_{ij})_{(i,j)=(1,1)}^{(k,\ell)}$, listing for each component $M_{e_j}(\CC)$ its multiplicity in $M_{d_i}(\CC)$. For $T$-linear derivations we can state the following result.
\begin{theorem}\label{DoubleDerivationsRelative}
Let $T=M_{e_1}(\CC)\oplus\dots M_{e_\ell}(\CC)$ be a finite dimensional semi simple subalgebra of $S$ with Bratelli diagram with respect to $S$ given by $(a_{ij})_{(i,j)=(1,1)}^{(k,\ell)}$, then
$$\DDer_T(S) \cong \bigoplus_{i=1}^k M_{d_i}(\CC)^{\oplus r_i} \oplus
\bigoplus_{i\neq j} M_{d_i\times d_j}(\CC)^{\oplus r_{ij}}$$
as $S$-modules, with $r_i = \sum_{u=1}^\ell a_{iu}^2 - 1$ and $r_{ij} = \sum_{u=1}^\ell a_{iu}a_{ju}$.
\end{theorem}
\begin{proof}
We use the identification from Theorem \ref{DoubleDerivationsNotRelative}. Write $f_{pq}^{rs}(i)$ for the basis elements in $M_{d_i}(\CC)^{d_i^2-1}$ ($1\leq r,s\leq d_i$, $1\leq p\leq d_i$, $2\leq q\leq d_i$) and $f_{pq}^{rs}(i,j)$ for the basis elements in $M_{d_i\times d_j}(\CC)^{d_id_j}$ ($i\neq j$, $1\leq p,s \leq d_j$, $1\leq q,r\leq d_i$).  Let $x\in\DDer(S)$. Let $y ^u_{vw}$ be the $(v,w)$-th basis element of the $u$-th component of $T$, then $y^u_{vw}$ is embedded in the $i$-th component of $S$ as
$$\sum_{h=0}^{a_{iu}-1} e_{(n_{iu}+he_u+v)(n_{iu}+he_u+w)},$$
where $n_{iu} = a_{i1}e_1+\dots a_{i(u-1)}e_{u-1}$ for $i>1$ and $n_{i1} = 0$. Then the term of $x(y^u_{vw})$ in $S_i\otimes S_j$ for $i\neq j$ equals
\begin{align}
&\sum_{r=1}^{n_{iu}}\sum_{(p,q)=(1,1)}^{(d_i,d_j)}\sum_{h=0}^{a_{ju}-1} x_{pq}^{r(n_{ju}+v+he_{u})}(i,j)f_{pq}^{r(n_{ju}+w+he_u)}(i,j)\label{Collection1}\\
+&\sum_{r=n_{iu}+1}^{n_{iu}+a_{iu}e_u}\sum_{(p,q)=(1,1)}^{(d_i,d_j)}\sum_{h=0}^{a_{ju}-1} x_{pq}^{r(n_{ju}+v+he_{u})}(i,j)f_{pq}^{r(n_{ju}+w+he_u)}(i,j)\label{Collection2}\\
-&\sum_{s=1}^{d_j}\sum_{(p,q)=(1,1)}^{(d_i,d_j)}\sum_{h=0}^{a_{iu}-1} x_{pq}^{(n_{iu}+v+he_u)s}(i,j)f^{(n_{iu}+w+he_u)s}_{pq}(i,j)\label{Collection3}\\
+&\sum_{r=n_{i(u+1)}+1}^{d_i}\sum_{(p,q)=(1,1)}^{(d_i,d_j)}\sum_{h=0}^{a_{ju}-1} x_{pq}^{r(n_{ju}+v+he_{u})}(i,j)f_{pq}^{r(n_{ju}+w+he_u)}(i,j).\label{Collection4}
\end{align}
Now $x\in \DDer_T(S)$ if and only if $x(y) = 0$ for all $y\in T$. Letting $u$, $v$ and $w$ run over all possible values, lines (\ref{Collection1}) and (\ref{Collection4}) in the expression above then yield that for $r\in [n_{iu}+1,n_{i(u+1)}]$ and $s\not\in [n_{ju}+1,n_{j(u+1)}]$ the term of $x$ in the $(r,s)$-th component of $M_{d_i\times d_j}(\CC)^{\oplus d_id_j}$ is zero. Choosing $v=w$ and running over all possible values, lines (\ref{Collection2}) and (\ref{Collection3}) imply that for $r\in [n_{iu}+1,n_{i(u+1)}]$ and $s\in [n_{ju}+1,n_{j(u+1)}]$ with $r=n_{iu}+xe_u + y$ and $s=n_{ju}+ae_u + b$, $y,b<e_u$ the term in the $(r,s)$-th component of $M_{d_i\times d_j}(\CC)^{\oplus d_id_j}$ is zero unless $y=b$. Choosing $v\neq w$ and letting these indices run over all possible values again then yields that for $y = b$ we have that the term in the $(r,s)$-th component of $M_{d_i\times d_j}(\CC)^{\oplus d_id_j}$ is equal to the term in the $(r',s)$-th component for $r' = n_{iu}=xe_u+y'$. This means that we only have $a_{iu}a_{ju}$ nonzero components of $x$ for each $u$, yielding the multiplicities $r_{ij}$ of the theorem.

The multiplicities $r_i$ are obtained through a completely analogous reasoning, taking into account that because of Theorem \ref{DoubleDerivationsNotRelative} we start with $d_i^2-1$ copies of $M_{d_i}(\CC)$ instead of $d_id_i$ copies. 
\end{proof}

\section{Double Poisson Tensors and Cycles}\label{DPTC}
Because a finite dimensional semi-simple algebra is formally smooth, all double Poisson brackets are determined by linear combinations of degree two elements $\delta\Delta$ in $\DDQ$ that are nilpotent (modulo commutators) with respect to the Schouten bracket on $\DDQ$ \cite[Proposition 4.1.2]{MichelDPA}. The double bracket corresponding to such a degree $2$ element being defined as
$$\db{a,b} := \Delta(b)'\delta(a)'' \otimes \delta(a)'\Delta(b)'' - \delta(b)'\Delta(a)''\otimes\Delta(a)'\delta(b)''.$$

Recall from \cite{MichelDPA} that the double Schouten bracket $\db{-,-}_s$ on $\DDQ$ is  defined as
\begin{eqnarray*}
\db{s,t}_s & = & 0 \\
\db{\delta,s}_s & =& \delta(s) \\
\db{\delta,\Delta}_s & = & \db{\delta,\Delta}_l + \db{\delta,\Delta}_r
\end{eqnarray*}
where
$$\db{\delta,\Delta}_l = \tau_{(23)}\db{\delta,\Delta}_{\tilde{l}} ~\mathrm{and}~
\db{\delta,\Delta}_r = \tau_{(12)}\db{\delta,\Delta}_{\tilde{r}}$$
with
$$\db{\delta,\Delta}_{\tilde{l}} := (\delta\otimes 1)\Delta - (1\otimes\Delta)\delta ~\mathrm{and}~
\db{\delta,\Delta}_{\tilde{r}} := (1\otimes\delta)\Delta - (\Delta\otimes 1)\delta$$
viewed as elements of $\DDer(S)\otimes S$ respectively $S\otimes\DDer(S)$. This bracket is a \emph{double Gerstenhaber bracket}; that is, it satisfies
\begin{enumerate}
\item a graded derivation property
$$\db{a,bc} = (-1)^{(| a| - 1)| b|}b\db{a,c} + \db{a,b}c;$$
\item a graded anti-symmetry property
$$\db{a,b} = -\sigma_{(12)}(-1)^{(|a|-1)(|b|-1)}\db{b,a},$$
where $$\sigma_s(a) := (-1)^t a_{s^{-1}(1)}\otimes\dots\otimes a_{s^{-1}(n)},$$
for $s\in\mathfrak{S}_n$ and $a = a_1\otimes\dots\otimes a_n$;
\item a graded double Jacobi identity
\begin{align*}
&\db{a_1,\db{a_2,a_3}'}\otimes \db{a_2,a_3}'' + (-1)^{(|a_1|-1)(| a_2| +| a_3|)}\sigma_{(123)}{\db{a_2,\db{a_3,a_1}'}\otimes\db{a_3,a_1}''} \\
&+ (-1)^{(| a_3|-1)(| a_1| +| a_2|)}\sigma_{(132)}{\db{a_3,\db{a_1,a_2}'}\otimes\db{a_1,a_2}''} = 0.
\end{align*}
\end{enumerate}

The Schouten bracket on $\DDQ$ becomes
\begin{prop}\label{SchoutenGenerators}
For the generators of $\DDQ$ defined in the previous section, we have
\begin{eqnarray*}
\db{x^i_{pq}, x^j_{rs}}_s 
& = & \left\{\begin{array}{rll} & 0 & i \neq j \\
&e^i_{p1}\otimes e_{1q}x_{rs}^i - x_{rs}^ie^i_{p1}\otimes e^i_{1q}& \\
+& e_{1s}^i\otimes x_{pq}^ie_{r1}^i - e_{1s}^ix_{pq}^i\otimes e_{r1}^i & i = j \\
 +& \delta_{rq}x_{ps}^i\otimes e_{11}^i - \delta_{ps}e_{11}^i\otimes x^i_{rq} & \end{array}\right .
\\
\db{x^i_{pq}, y^{uv}_{rs}}_s 
& = & \left\{\begin{array}{ll} 0 & i \neq u,v \\
e_{1s}^v\otimes x_{pq}^ie_{r1}^i + \delta_{qr}y_{ps}^{iv}\otimes e_{11}^i - y^{iv}_{rs}e_{p1}^i\otimes e_{1q}^i & u = i \\
e_{p1}^i\otimes e_{1q}^iy_{rs}^{ui} - \delta_{ps}e_{11}^i\otimes y_{rq}^{ui} - e_{1s}^ix_{pq}^i\otimes e_{r1}^u
& v = i
\end{array}\right .
\\ 
\db{y_{pq}^{rs},y_{ab}^{cd}}_s & = & 
\left\{
\begin{array}{ll}
e_{p1}^r\otimes e_{1q}^sy_{ab}^{cd} + e_{1b}^d\otimes y_{pq}^{rs}e_{a1}^c -e_{1b}^dy_{pq}^{rs}\otimes e_{a1}^c - y_{ab}^{cd} e_{p1}^r\otimes e_{1q}^s & r\neq d, c\neq s \\
- \delta_{bp}e_{11}^d\otimes y_{aq}^{cs} & r = d, c\neq s \\ 
- \delta_{bp}e_{11}^r\otimes x_{aq}^{s} + \delta_{qa}x_{pb}^{r}\otimes e_{11}^s &  r =  d, c =  s
\end{array}\right.
\end{eqnarray*}
where we use the shorthand notation $x_{11}^i = -\sum_{r= 2}^{d_i} x_{rr}^i$.
\end{prop}
\begin{proof}
We show the first equality holds. The other computations are analogous. By definition of the double Schouten bracket it is obvious that $\db{x^i_{pq}, x^j_{rs}} = 0$ if $i\neq j$. Now assume $i = j$ and denote $x^i_{pq} = x_{pq}$ and $x^i_{rs} = x_{rs}$, then for $z\in M_{d_i}(\CC)$ we get
\begin{eqnarray*}
(1\otimes x_{pq})x_{rs}(z) - (x_{rs}\otimes 1)x_{pq}(z) & = & 
e_{r1}\otimes e_{p1}\otimes e_{1q}e_{1s}z - e_{r1}\otimes e_{1s}ze_{p1}\otimes e_{1q} \\
& & -ze_{r1}\otimes e_{p1}\otimes e_{1q}e_{1s} + ze_{r1}\otimes e_{1s}e_{p1}\otimes e_{1q} \\
& & -e_{r1}\otimes e_{1s}e_{p1}\otimes e_{1q}z + e_{p1}e_{r1}\otimes e_{1s}\otimes e_{1q}z \\
& & +e_{r1}\otimes e_{1s}ze_{p1}\otimes e_{1q} - ze_{p1}e_{r1}\otimes e_{1s}\otimes e_{1q}
\end{eqnarray*}
Applying $\tau_{(12)}$ to this expression yields
$$e_{p1}\otimes e_{1q}x_{rs}(z) + e_{1s}\otimes x_{pq}(z)e_{r1} - \delta_{ps}e_{11}\otimes x_{rq}(z).$$
Changing the sign, interchanging the indices and applying to the previous computation $\tau_{(23)}$ then yields
$$-e_{1s}x_{pq}(z)\otimes _{r1} - \delta_{rq}x_{ps}(z)\otimes e_{11} + x_{rs}(z)e_{p1}\otimes e_{1q}.$$
Adding this expression to the one found in the previous paragraph then yields the expression in the statement of the proposition.
\end{proof}
For the remainder of the paper, the index $(1,1)$ at a loop will be the shorthand notation introduced in the preceding proposition.

For $S = \CC^{\oplus n}$ this structure descends to the following double Gerstenhaber bracket on the quiver $Q_S$.
\begin{cor}
The double Gerstenhaber bracket on $\CC Q_S$ for $S = \CC^{\oplus n}$ is determined by
\begin{eqnarray*}
\db{i\leftarrow j, k \leftarrow i} & = & -i\otimes (k\leftarrow j) ~(j\neq k)\\
\db{k \leftarrow i,i\leftarrow j} & = & (k\leftarrow j)\otimes i~(j\neq k)\\
\db{i\leftarrow j, k \leftarrow \ell} & = & 0 ~(\textrm{otherwise})
\end{eqnarray*}
\end{cor}
More generally, we know $\DDQ/[\DDQ,\DDQ][1]$ is a graded Lie algebra with graded Lie bracket $\{-,-\}_s$. This algebra can be described in terms of the double derivation quiver, where the path algebra of this quiver is considered graded with the arrows of degree $1$ and the orthogonal idempotents of degree $0$.

\begin{theorem}\label{DDQGradedLie}
For a semi-simple algebra $S$, we have that $\DDQ/[\DDQ,\DDQ][1]$ with bracket $\{-,-\}_s$ is isomorphic as a graded Lie algebra to $\CC Q_S/[\CC Q_S,\CC Q_S]_{super}$ with bracket defined as follows
$$
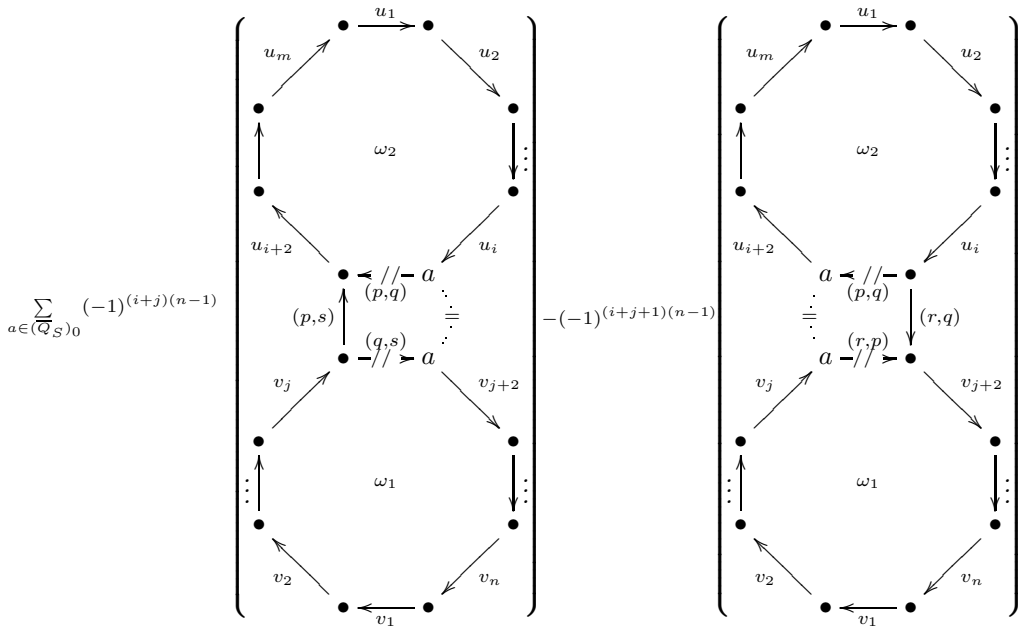
\xymatrix@R=1.75pc@C=1.75pc{
\save[0,1].[7,1]!C*+<10pt,0pt>\frm{(}\restore
\save[0,1].[7,1]!C+<-55pt,0pt>
\drop{\scriptstyle\sum\limits_{a\in (\overline{Q}_S)_0} (-1)^{(i+j)(n-1)}}
\restore
\save[0,4].[7,4]*+<10pt,0pt>\frm{)}\restore
\save[0,7].[7,7]*+<10pt,0pt>\frm{(}\restore
\save[0,7].[7,7]!C+<-42pt,0pt>
\drop{\scriptstyle-(-1)^{(i+j+1)(n-1)}}
\restore
\save[0,10].[7,10]*+<10pt,0pt>\frm{)}\restore
&
&\bullet \ar[r]^{u_1} & \bullet \ar[dr]^{u_2} 
&&& 
& &\bullet \ar[r]^{u_1} & \bullet \ar[dr]^{u_2} & 
\\ & 
\bullet\ar[ur]^{u_m}&& & \bullet \ar[d]^{\vdots}
&&&
\bullet\ar[ur]^{u_m}&& & \bullet \ar[d]^{\vdots}
\\ &
\bullet\ar[u]\ar@{}[urrr]|{\omega_2}& & & \bullet \ar[dl]^{u_i}
&&&
\bullet\ar[u]\ar@{}[urrr]|{\omega_2}&& & \bullet \ar[dl]^{u_i}
\\ &
 &\bullet\ar[ul]^{u_{i+2}}& a \ar[l]^{(p,q)}|{~~//}\ar@/^.75pc/@{.}[d]|=&
&&&
&a \ar[ul]^{u_{i+2}} & \bullet\ar[l]^{(p,q)}|{~~//}\ar[d]^{(r,q)} &
 \\ &
 &\bullet \ar[u]^{(p,s)}\ar[r]^{(q,s)}|{//~~}  & a\ar[dr]^{v_{j+2}}& &
&& 
& a\ar[r]^{(r,p)}|{//~~}\ar@/^.5pc/@{.}[u]|{=} & \bullet\ar[dr]^{v_{j+2}}& 
\\   &
\bullet\ar[ur]^{v_j}&& & \bullet \ar[d]^{\vdots}
&& & 
\bullet\ar[ur]^{v_j}&& & \bullet \ar[d]^{\vdots}
\\ &
\bullet\ar[u]^{\vdots}\ar@{}[urrr]|{\omega_1}&& & \bullet \ar[dl]^{v_n}
&&& 
\bullet\ar[u]^{\vdots}\ar@{}[urrr]|{\omega_1}&& & \bullet \ar[dl]^{v_n}
\\ &
&\bullet\ar[ul]^{v_{2}}& \bullet \ar[l]^{v_1}
&&& 
& &\bullet\ar[ul]^{v_{2}}& \bullet \ar[l]^{v_1} & & 
\\
}
$$
That is, we fix representatives for $\omega_1$ and $\omega_2$ in $\CC Q_S$ and for every vertex $a\in (Q_S)_0$ we look for an occurrence of $a$ in $\omega_1$ and in $\omega_2$. If $a = h(u_i)$ in $\omega_1$ and $a = t(v_{j+2})$ in $\omega_2$ and the secondary colour of $u_{i+1}$ equals the primary colour of $v_{j+1}$, we remove $u_{i+1}$ from $\omega_1$ and $v_{j+1}$ from $\omega_2$, glue $\omega_1$ to $\omega_2$ in $a$ and connect the loose ends $t(v_{j+1})$ and $h(u_{i+1})$ with an arrow with primary colour equal to that of $u_{i+1}$ and secondary colour equal to that of $v_{j+1}$. The necklace thus obtained is multiplied by a factor $(-1)^{(i+j)(n-1)}$. An analogous procedure is followed each time $a$ is the tail of an arrow in $\omega_1$ and the head of an arrow in $\omega_2$, with the roles of primary and secondary colours interchanged and factor $-(-1)^{(i+j+1)(n-1)}$.
\end{theorem}
\begin{proof}
First of all, note that $\DDQ/[\DDQ,\DDQ]$ is generated by all cycles for $Q_S$ because for generators $u_i \in \{x_{pq}^j,y_{rs}^{k\ell}\}$ and elements $c_i \in S$ we have for any word
\begin{eqnarray*}
c_1u_1c_2u_2\dots c_pu_pc_{p+1} & = & c_1e_{11}^{h(u_1)}u_1 e_{11}^{t(u_1)}c_2\dots c_1e_{11}^{h(u_p)}u_p e_{11}^{t(u_p)}c_{p+1} \\
& = & (c_1)_{11}^{h(u_1)}\dots (c_p)_{11}^{h(u_p)} (c_{p+1})_{11}^{h(u_1)}u_1\dots u_p \mathrm{~mod~} [\DDQ,\DDQ].
\end{eqnarray*}
Now let $u_1u_2\dots u_m$ and $v_1v_2\dots v_n$ be two necklaces in $\DDQ$, then
\begin{align*}
&\db{v_1v_2\dots v_n,u_1u_2\dots u_m} =\\
& \sum_{i=0}^{m-1}\sum_{j = 0}^{n-1} -(-1)^{i(n-1)} u_1\dots u_i \sigma_{(12)}(v_1\dots v_j\db{u_{i+1},v_{j+1}}v_{j+2}\dots v_n)u_{i+2}\dots u_m.
\end{align*}
We have four different settings to consider. First of all assume $u_{i+1} = x_{pq}^a$ and $v_{j+1} = x_{rs}^b$. Proposition \ref{SchoutenGenerators} then yields $\db{u_{i+1},v_{j+1}} = 0$ if $a\neq b$ and
\begin{align*}
&  u_1\dots u_i \sigma_{(12)}(v_1\dots v_j\db{u_{i+1},v_{j+1}}v_{j+2}\dots v_n)u_{i+2}\dots u_m = \\
& (-1)^{j(n-j)}  u_1\dots u_i e_{1q}x_{rs}^av_{j+2}\dots v_n\otimes v_1\dots v_je^a_{p1}u_{i+2}\dots u_m \\
& - (-1)^{(j+1)(n-j-1)}  u_1\dots u_i e^a_{1q}v_{j+2}\dots v_n\otimes v_1\dots v_jx_{rs}^ae^a_{p1}u_{i+2}\dots u_m \\
&+(-1)^{j(n-j)} u_1\dots u_i x_{pq}^ae_{r1}^av_{j+2}\dots v_n\otimes v_1\dots v_je_{1s}^iu_{i+2}\dots u_m \\
&- (-1)^{(j+1)(n-j-1)} u_1\dots u_i e_{r1}^av_{j+2}\dots v_n\otimes v_1\dots v_je_{1s}^ax_{pq}^au_{i+2}\dots u_m \\
&+ (-1)^{(j+1)(n-j-1)}\delta_{rq} u_1\dots u_i e_{11}^av_{j+2}\dots v_n\otimes v_1\dots v_jx_{ps}^au_{i+2}\dots u_m \\
&-(-1)^{j(n-j)}\delta_{ps} u_1\dots u_i x^a_{rq}v_{j+2}\dots v_n\otimes v_1\dots v_je_{11}^au_{i+2}\dots u_m.
\end{align*}
Now by definition $(p,q)\neq (1,1)$ and $(r,s)\neq (1,1)$, so modulo commutators this expression is mapped by the multiplication to
\begin{align*}
&\delta_{rq}(-1)^{(j+1)(n-j-1)}u_1\dots u_iv_{j+2}\dots v_nv_1\dots v_jx_{ps}^au_{i+2}\dots u_m \\
&-\delta_{ps} (-1)^{j(n-j)}u_1\dots u_i x^a_{rq}v_{j+2}\dots v_nv_1\dots v_ju_{i+2}\dots u_m.
\end{align*}
Next, assume $u_{i+1} = x_{pq}^a$ and $v_{i+1} = y_{rs}^{bc}$. This yields $\db{u_{i+1},v_{j+1}} = 0$ if $a\neq b,c$ and a similar computation as the previous yields
\begin{align*}
&  u_1\dots u_i \sigma_{(12)}(v_1\dots v_j\db{u_{i+1},v_{j+1}}v_{j+2}\dots v_n)u_{i+2}\dots u_m \mapsto \\
&\delta_{qr}(-1)^{(j+1)(n-j-1)}u_1\dots u_i v_{j+2}\dots v_nv_1\dots v_j y_{ps}^{bc}u_{i+2}\dots u_m
\end{align*}
if $a = b$ and
\begin{align*}
&  u_1\dots u_i \sigma_{(12)}(v_1\dots v_j\db{u_{i+1},v_{j+1}}v_{j+2}\dots v_n)u_{i+2}\dots u_m \mapsto \\
&-\delta_{ps}(-1)^{j(n-j)}u_1\dots u_i v_{j+2}\dots v_nv_1\dots v_j y_{rq}^{bc}u_{i+2}\dots u_m
\end{align*}
if $a = c$. For  $u_{i+1} = y_{rs}^{bc}$ and $v_{i+1} = x_{pq}^a$ we get
\begin{align*}
&  u_1\dots u_i \sigma_{(12)}(v_1\dots v_j\db{u_{i+1},v_{j+1}}v_{j+2}\dots v_n)u_{i+2}\dots u_m \mapsto \\
&-\delta_{qr}(-1)^{j(n-j)}u_1\dots u_i v_{j+2}\dots v_nv_1\dots v_j y_{ps}^{bc}u_{i+2}\dots u_m
\end{align*}
if $a = b$ and
\begin{align*}
&  u_1\dots u_i \sigma_{(12)}(v_1\dots v_j\db{u_{i+1},v_{j+1}}v_{j+2}\dots v_n)u_{i+2}\dots u_m \mapsto \\
&\delta_{ps}(-1)^{(j+1)(n-j-1)}u_1\dots u_i v_{j+2}\dots v_nv_1\dots v_j y_{rq}^{bc}u_{i+2}\dots u_m
\end{align*}
if $a = c$. Finally, for $u_{i+1} = y_{pq}^{ab}$ and $v_{i+1} = y_{rs}^{cd}$ we get $\{u_{i+1},v_{i+1}\}_s = 0$ unless $b = c$ or $a = d$ in which case we get
either
\begin{align*}
&  u_1\dots u_i \sigma_{(12)}(v_1\dots v_j\db{u_{i+1},v_{j+1}}v_{j+2}\dots v_n)u_{i+2}\dots u_m \mapsto \\
&- \delta_{qr}(-1)^{j(n-j)}u_1\dots u_iy_{ps}^{cb}v_{j+2}\dots v_nv_1\dots v_ju_{i+2}\dots u_m,
\end{align*}
if $a = d$ and $b\neq c$,
\begin{align*}
&  u_1\dots u_i \sigma_{(12)}(v_1\dots v_j\db{u_{i+1},v_{j+1}}v_{j+2}\dots v_n)u_{i+2}\dots u_m \mapsto \\
&\delta_{ps}(-1)^{(j+1)(n-j-1)}u_1\dots u_iv_{j+2}\dots v_nv_1\dots v_jy_{rq}^{ad}u_{i+2}\dots u_m,
\end{align*}
if $b = c$ and $a\neq d$ or
\begin{align*}
&  u_1\dots u_i \sigma_{(12)}(v_1\dots v_j\db{u_{i+1},v_{j+1}}v_{j+2}\dots v_n)u_{i+2}\dots u_m \mapsto \\
&-\delta_{sp}(-1)^{j(n-j)}u_1\dots u_ix_{rq}^bv_{j+2}\dots v_nv_1\dots v_j u_{i+2}\dots u_m \\
&+\delta_{qr}(-1)^{(j+1)(n-j-1)}u_1\dots u_iv_{j+2}\dots v_n v_1\dots v_ix^a_{ps} u_{i+2}\dots u_m.
\end{align*}
if $b = c$ and $a = d$.
\end{proof}
Now, in order to determine which length two elements in $\DDQ$ yield nontrivial double Poisson brackets on $S$, note that
\begin{lemma}
For $P = x^i_{pq}y^{rs}_{uv}$ or $P = y^{rs}_{uv}x^i_{pq}$ with $i,p,q,r,s,u,v$ arbitrary and for $P= y_{pq}^{rs}y_{tu}^{vw}$ with $p,q,t,u$ arbitrary and $s\neq v$ or $w\neq r$ we have that
$$\db{-,-}_P = 0.$$
\end{lemma}
\begin{proof}
This is due to the fact that in the definition of $\db{-,-}_P$ the tensor products are obtained through componentwise multiplication.
\end{proof}
This lemma, in combination with the proposition preceding it, proves that
\begin{prop}\label{OnlyCycles}
Let $\db{-,-}$ be a nonzero double bracket on $S$, then $\db{-,-}$ is completely determined by a linear combination
$$\sum \alpha_{abcd}^{pq}y_{ab}^{pq}y_{cd}^{qp} + \beta_{efgh}^ix_{ef}^ix_{gh}^i$$
with all $\alpha_{abcd}^{pq}$ and $\beta_{efgh}^i$ in $\CC$.
\end{prop}

The next two lemmas determine which cycles of length two yield double Poisson brackets. 
\begin{lemma}\label{DoubleDerivationCycle1}
For $P = y_{ab}^{pq}y_{cd}^{qp}$ we have that
\begin{enumerate}
\item if $S_p\neq\CC$ and $S_q\neq \CC$ then
$$\{P,P\} = 0 ~\textrm{mod}~ [\DDQ,\DDQ] \Leftrightarrow (a-d)(b-c)\neq 0 $$
\item if $S_p\neq \CC$ and $S_q = \CC$ then
$$\{P,P\} = 0 ~\textrm{mod}~ [\DDQ,\DDQ] \Leftrightarrow a-d\neq 0$$
\item if $S_p = \CC$ and $S_q \neq \CC$ then
$$\{P,P\} = 0 ~\textrm{mod}~ [\DDQ,\DDQ] \Leftrightarrow b-c\neq 0$$
\item if $S_p = \CC$ and $S_q = \CC$  then
$$\{P,P\} = 0 ~\textrm{mod}~ [\DDQ,\DDQ]$$
\end{enumerate}
and in each of these cases $\db{-,-}_P$ is nontrivial.
\end{lemma}
\begin{proof}
Using the description of the bracket on $\CC\overline{Q}_S/[\CC\overline{Q}_S,\CC\overline{Q}_S]$ from Proposition \ref{DDQGradedLie}, we compute
$$\left\{
\xymatrix@R=1.75pc@C=1.75pc{
p\ar@/^.5pc/[r]^{(c,d)} & q\ar@/^.5pc/[l]^{(a,b)}
},
\xymatrix@R=1.75pc@C=1.75pc{
p\ar@/^.5pc/[r]^{(c,d)} & q\ar@/^.5pc/[l]^{(a,b)}
}
\right\} = 
2\delta_{ad}\vcenter{\xymatrix@R=.90pc@C=.90pc{
p\ar[rr]^{(a,b)} & & q\ar[ddl]^{(c,d)} \\ \\ & p\ar[uul]^{(c,b)}
} }
+ 2\delta_{bc}\vcenter{\xymatrix@R=.90pc@C=.90pc{
q\ar[rr]^{(c,d)} & & p\ar[ddl]^{(a,b)} \\ \\ & q\ar[uul]^{(a,d)}
}}.$$
Which is nonzero if $a=d$ and the loop at $p$ exists or if $b=c$ and the loop at $p$ exists, yielding the four situations described in the lemma. The fact that the bracket obtained from this double Poisson tensor is non-zero is easy to verify.
\end{proof}

\begin{lemma}\label{DoubleDerivationCycle2}
For $P = x^u_{pq}x^u_{rs}$ we have that
$$\{P,P\} = 0 ~\textrm{mod}~ [\DDQ,\DDQ]$$
and $\db{-,-}_P$ is nontrivial. If and only if either $(p-q)(p-s)(r-s)(r-q)\neq 0$, or  $p = q = r$, or $r = s = p$.
\end{lemma}
\begin{proof}
Again we use the description of the bracket from Proposition \ref{DDQGradedLie}.
We have
\begin{eqnarray*}
\left\{
\xymatrix{
u\ar@/^.5pc/[r]^{(p,q)} & u\ar@/^.5pc/[l]^{(r,s)}
},
\xymatrix{
u\ar@/^.5pc/[r]^{(p,q)} & u\ar@/^.5pc/[l]^{(r,s)}
}
\right\}
& = & 
2\delta_{ps}\vcenter{\xymatrix@R=.90pc@C=.90pc{
u \ar[rr]^{(p,q)} & & u\ar[ddl]^{(r,s)} \\ \\ & u\ar[uul]^{(r,q)}
}}
-2\delta_{rs}\vcenter{\xymatrix@R=.90pc@C=.90pc{
u \ar[rr]^{(p,q)} & & u\ar[ddl]^{(p,q)} \\ \\ & u\ar[uul]^{(r,r)}
}}
\\ & & 
+2\delta_{rq}\vcenter{\xymatrix@R=.90pc@C=.90pc{
u \ar[rr]^{(p,q)} & & u\ar[ddl]^{(p,s)} \\ \\ & u\ar[uul]^{(r,s)}
}}
-2\delta_{pq}\vcenter{\xymatrix@R=.90pc@C=.90pc{
u \ar[rr]^{(p,p)} & & u\ar[ddl]^{(r,s)} \\ \\ & u\ar[uul]^{(r,s)}
}}
\end{eqnarray*}
Which is exactly nonzero when the conditions formulated in the lemma are satisfied. Again, nontriviality is straightforward.
\end{proof}

For $S = \CC^{\oplus n}$, we can easily show the following
\begin{prop}\label{NoJointVertices}
Let 
$$P = \sum_{i<j} \alpha_{ij}y^{ij}y^{ji} \in \DDQ,$$
then $P$ determines a double Poisson bracket on $S$ if and only if the following relation between the $\alpha$ holds:
$$\forall i<j<k: \alpha_{ij}\alpha_{ik} + \alpha_{ik}\alpha_{jk} - \alpha_{ij}\alpha_{jk} = 0.$$
\end{prop}
\begin{proof}
A straightforward computation shows that for $i,j,k,\ell$ pairwise different we have
\begin{eqnarray*}
\{y^{ij}y^{ji},y^{ki}y^{ik}\} & = & -y^{ki}y^{jk}y^{ij} - y^{ji}y^{kj}y^{ik} ~\textrm{mod}~ [\DDQ,\DDQ]\\
\{y^{ij}y^{ji},y^{jk}y^{kj}\} & = & -y^{jk}y^{ij}y^{ki} - y^{ik}y^{ji}y^{kj} ~\textrm{mod}~ [\DDQ,\DDQ]\\
\{y^{ij}y^{ji},y^{ik}y^{ki}\} & = & y^{ik}y^{ji}y^{kj} + y^{jk}y^{ij}y^{ki} ~\textrm{mod}~ [\DDQ,\DDQ]\\
\{y^{ji}y^{ij},y^{ki}y^{ik}\} & = & y^{ik}y^{ji}y^{kj} + y^{jk}y^{ij}y^{ki} ~\textrm{mod}~ [\DDQ,\DDQ]\\
\{y^{ij}y^{ji},y^{k\ell}y^{\ell k}\} & = & 0 ~\textrm{mod}~ [\DDQ,\DDQ]
\end{eqnarray*}
This yields
\begin{eqnarray*}
\{P,P\} & = & \sum_{i<j}\sum_{k<\ell} \alpha_{ij}\alpha_{k\ell}\{y^{ij}y^{ji},y^{k\ell}y^{\ell k}\} \\
& = & \sum_{i,k<j} \alpha_{ij}\alpha_{kj}(y^{jk}y^{ij}y^{ki} + y^{ik}y^{ji}y^{kj}) \\
&   & - \sum_{i<j<\ell} \alpha_{ij}\alpha_{j\ell}(y^{j\ell}y^{ij}y^{\ell i} + y^{ji}y^{\ell j}y^{i\ell}) \\
&   & + \sum_{i<j,\ell} \alpha_{ij}\alpha_{i\ell}( y^{i\ell}y^{ji}y^{\ell j} + y^{j\ell}y^{ij}y^{\ell i}) \\
&   &  - \sum_{k<i<j} \alpha_{ij}\alpha_{ki}(y^{ki}y^{jk}y^{ij} + y^{ji}y^{kj}y^{ik}) ~\textrm{mod}~ [\DDQ,\DDQ] \\
& = & \sum_{i<j<k} 2(\alpha_{ij}\alpha_{ik} + \alpha_{ik}\alpha_{jk} - \alpha_{ij}\alpha_{jk})(y^{ji}y^{kj}y^{ik} + y^{ki}y^{jk}y^{ij})~\textrm{mod}~ [\DDQ,\DDQ] .
\end{eqnarray*}
So we have a double Poisson bracket if and only if for all $i<j<k$ we have 
$$\alpha_{ij}\alpha_{ik} + \alpha_{ik}\alpha_{jk} - \alpha_{ij}\alpha_{jk} = 0.$$
\end{proof}

The cycles of length two can be seen as being analogous to the classical Poisson tensors $\frac{\partial}{\partial x_i}\wedge\frac{\partial}{\partial x_j}$. In classical Poisson geometry, these Poisson tensors play an important role in the Splitting Theorem \cite{Weinstein}, where any Poisson structure can be split into a symplectic part and a totally degenerate part, the totally degenerate part being of the form $\varphi_{ij}(x)\frac{\partial}{\partial x_i}\wedge\frac{\partial}{\partial x_j}$. If $\varphi_{ij}(x) = c_{ij}$, the totally degenerate part gives rise to an affine Poisson structure for any choice of $c_{ij}$ \cite{AffinePoisson}. In the case of double Poisson structures, the previous proposition shows that this no longer holds.

Another important concept introduced in \cite{MichelDPA} is that of a moment map for a double Poisson bracket. In the case of a finite dimensional semi-simple algebra $S = \CC^{\oplus n}$ we have
\begin{prop}\label{MomentMaps}
Let $S = \CC^{\oplus{n}}$ with orthogonal idempotents $e_i$ and let $P = \sum_{i<j} c_{ij}d_{ij}d_{ji}$ determine a double Poisson bracket with $d_{ij}$ the derivation corresponding to the unique arrow from $j$ to $i$, then there exists a moment map $\mu$ for $P$ if and only if all $c_{ij}$ are nonzero.
In this case, the moment map $\mu$ is unique up to a constant term and can be written as
$$\mu = -\sum_{i=2}^n \frac{1}{c_{1i}}e_i.$$
\end{prop}
\begin{proof}
First of all note that for the double Schouten bracket we have
\begin{eqnarray*}
\db{\delta\Delta,\mu} & = & \db{\mu,\delta\Delta}^{op} \\
& = & (-\delta\db{\mu,\Delta} + \db{\mu,\delta}\Delta)^{op} \\
&=& -\delta*\Delta(\mu) + \delta(\mu)*\Delta.
\end{eqnarray*}
So for $\mu = \sum_i \mu_ie_i$ we get
\begin{eqnarray*}
\db{P,\mu} & = & \sum_{i<j} c_{ij}\db{d_{ij}d_{ji},\mu} \\
&=& \sum_{i<j} c_{ij}(-d_{ij}*d_{ji}(\mu) + d_{ij}(\mu)d_{ji} \\
&=& \sum_{i<j} c_{ij}(-(\mu_i-\mu_j)e_j\otimes d_{ij} + (\mu_j-\mu_i)d_{ji}\otimes e_i \\
&=& \sum_{i<j} c_{ij}(\mu_j-\mu_i)(e_j\otimes d_{ij} + d_{ji}\otimes e_i)
\end{eqnarray*}
under the multiplication this maps to
$$\sum_{i<j} c_{ij}(\mu_j-\mu_i)(d_{ij}+d_{ji}),$$
which on $e_i$ is equal to
$$\sum_{i<j} c_{ij}(\mu_j-\mu_i)(-e_i\otimes e_j + e_j\otimes e_i) + \sum_{i>j} c_{ji}(\mu_j-\mu_i)(e_i\otimes e_j - e_j\otimes e_i).$$
Now 
$$E(e_i) = 1\otimes e_i - e_i\otimes 1 = \sum_{j\neq i} e_j\otimes e_i - e_i\otimes e_j,$$
so $\{P,\mu\} = - E$ implies $c_{ij}(\mu_j - \mu_i) = -1$ if $j>i$ and $c_{ij}(\mu_j - \mu_i) = 1$ if $j<i$. This means all $c_{ij}$ must be nonzero and for $i<j$ we have
$$\mu_i - \mu_j = \frac{1}{c_{ij}},$$
which means the $c_{ij}$ have to satisfy the additional relation for $i<j<k$:
$$c_{ij}c_{jk} = c_{ij}c_{jk} - c_{ik}c_{jk}.$$
But as $P$ determines a double Poisson bracket, by Proposition \ref{NoJointVertices} we know this condition is automatically satisfied.
The uniqueness follows from the fact that $e_1+\dots+e_n = 1$ and the fact that one $\mu_i$ determines all other coefficients.
\end{proof}

\section{Double Poisson-Lichnerowicz Cohomology}\label{DPLC}
In \cite{Lich}, Lichnerowicz observed that $d_\pi = \{\pi,-\}$ with $\pi$ a Poisson tensor for a Poisson manifold $M$ is a square $0$ derivation of degree $+1$, which yields a complex
$$0\stackrel{d_\pi}{\rightarrow} \mathcal{O}(M) \stackrel{d_\pi}{\rightarrow} D\mathcal{O}(M) \stackrel{d_\pi}{\rightarrow} \wedge^2 D\mathcal{O}(M)\stackrel{d_\pi}{\rightarrow} \dots,$$
the homology of which is called the \emph{Poisson-Lichnerowicz cohomology}. In this section, we show there is an analogous cohomology on $\DDQ$ that descends to the classical Poisson-Lichnerowicz cohomology on the quotient spaces of the representation spaces of the algebra.

Let $A$ be an associative algebra with unit. From \cite[\S 7]{MichelDPA} we know that the Poisson bracket on $\rep_{n}(A)$ and $\iss_n(A)$ induced by a double Poisson tensor $P$ corresponds to the Poisson tensor $tr(P)$. We furthermore know that the map $tr:\DDQ/[\DDQ,\DDQ][1]\rightarrow\bigwedge\Der(\mathcal{O}(\rep_n(A))$ is a morphism of graded Lie algebras, so we have a morphism of complexes
$$\xymatrix{
0 \ar[r] & (\DDQ/[\DDQ,\DDQ])^0\ar[d]^{tr}\ar[r]^{d_P} &(\DDQ/[\DDQ,\DDQ])^1\ar[d]^{tr}\ar[r]^{d_P} &(\DDQ/[\DDQ,\DDQ])^2\ar[d]^{tr}\ar[r] & \dots \\
0 \ar[r] & \mathcal{O}(\rep_n(A)) \ar[r]^{d_{tr(P)}} & \Der(\mathcal{O}(\rep_n(A))) \ar[r]^{d_{tr(P)}} & \wedge^2\Der(\mathcal{O}(\rep_n(A))) \ar[r] & \dots
}$$
which restricts to a morphism of complexes
$$\xymatrix{
0 \ar[r] & (\DDQ/[\DDQ,\DDQ])^0\ar[d]^{tr}\ar[r]^{d_P} &(\DDQ/[\DDQ,\DDQ])^1\ar[d]^{tr}\ar[r]^{d_P} &(\DDQ/[\DDQ,\DDQ])^2\ar[d]^{tr}\ar[r] & \dots \\
0 \ar[r] & \mathcal{O}(\iss_n(A)) \ar[r]^{d_{tr(P)}} & \Der(\mathcal{O}(\iss_n(A))) \ar[r]^{d_{tr(P)}} & \wedge^2\Der(\mathcal{O}(\iss_n(A))) \ar[r] & \dots
}$$
So there is a map from the homology $H^\bullet_P(A)$ of the upper chain complex, which we call the \emph{double Poisson-Lichnerowicz cohomology} to the classical Poisson-Lichnerowicz cohomology on $\rep_n(A)$ and $\iss_n(A)$.

Using the description of $\DDQ/[\DDQ,\DDQ][1]$ from Proposition \ref{DDQGradedLie} we can compute the double Poisson-Lichnerowicz cohomology for a semi-simple algebra $S$ when the double Poisson bracket is given by a single necklace in a straightforward way. We illustrate this for the zero-th and first cohomology groups.  It is easy to see that $H^0_P(S) = 0$ from the description of the double Schouten bracket. For $H^1_P(S)$ we get
\begin{prop}\label{DPLH1}
For a double Poisson bracket $P = \xymatrix@R=1.5pc@C=1.5pc{i\ar@/^.5pc/[r]^{(p,q)} & j\ar@/^.5pc/[l]^{(r,s)}}$ with $q\neq r$ and $p\neq s$ we have that $H^1_P(S)$ is generated by
\begin{itemize}
\item all loops of all bi-colours in vertices $k\neq i,j $ in $\overline{Q}_S$,
\item the loops in $i$ with all possible bi-colourings where the primary colour is different from $q$ and the secondary colour is different from $r$,
\item the loops in $j$ with all possible bi-colourings where the primary colour is different from $s$ and the secondary colour is different from $p$,
\item the sum $\xymatrix{i\ar@(ul,dl)_{(r,r)}} + \xymatrix{j\ar@(ur,dr)^{(s,s)}}$ and
\item the sum $\xymatrix{i\ar@(ul,dl)_{(p,p)}} + \xymatrix{j\ar@(ur,dr)^{(q,q)}}$,
\end{itemize}
and hence has dimension $\sum\limits_{k\neq i,j} (n_k^2-1) + (n_i-1)^2 + (n_j-1)^2$.

For a double Poisson bracket $P = \xymatrix@R=1.5pc@C=1.5pc{i\ar@/^.5pc/[r]^{(p,q)} & i\ar@/^.5pc/[l]^{(r,s)}}$ we have that $H^1_P(S)$ is generated by
all loops of all bi-colours in vertices $j\neq i$ and the loop in $i$ with all possible bi-colours where the primary colour differs from $s$ and $q$  and the secondary colour differs from $p$ and $r$ if $(p-q)(p-s)(r-s)(r-q)\neq 0$.
\end{prop}
\begin{proof}
The first claim follows immediately from the fact that
\begin{eqnarray*}
\left\{
\xymatrix@R=1.5pc@C=1.5pc{i\ar@/^.5pc/[r]^{(p,q)} & j\ar@/^.5pc/[l]^{(r,s)}},
\xymatrix@R=1.5pc@C=1.5pc{i\ar[r]^{(u,v)} & i}
\right\}
& = & 
\delta_{vr}\xymatrix@R=1.5pc@C=1.5pc{i\ar@/^.5pc/[r]^{(p,q)} & j\ar@/^.5pc/[l]^{(u,s)}}
-
\delta_{uq}\xymatrix@R=1.5pc@C=1.5pc{i\ar@/^.5pc/[r]^{(p,v)} & j\ar@/^.5pc/[l]^{(r,s)}},
\end{eqnarray*}
\begin{eqnarray*}
\left\{
\xymatrix@R=1.5pc@C=1.5pc{i\ar@/^.5pc/[r]^{(p,q)} & j\ar@/^.5pc/[l]^{(r,s)}},
\xymatrix@R=1.5pc@C=1.5pc{j\ar[r]^{(u',v')} & j}
\right\}
& = & 
-\delta_{su'}\xymatrix@R=1.5pc@C=1.5pc{i\ar@/^.5pc/[r]^{(p,q)} & j\ar@/^.5pc/[l]^{(r,v')}}
+
\delta_{v'p}\xymatrix@R=1.5pc@C=1.5pc{i\ar@/^.5pc/[r]^{(u',q)} & j\ar@/^.5pc/[l]^{(r,s)}},
\end{eqnarray*}
and
\begin{eqnarray*}
\left\{
\xymatrix@R=1.5pc@C=1.5pc{i\ar@/^.5pc/[r]^{(p,q)} & j\ar@/^.5pc/[l]^{(r,s)}},
\xymatrix@R=1.5pc@C=1.5pc{k\ar[r]^{(u,v)} & k}
\right\}
& = & 0
\end{eqnarray*}
if $k \neq i,j$, in combination with the fact that $\{P, S\} = 0$. The second claim follows immediately from
\begin{eqnarray*}
\left\{
\xymatrix@R=1.5pc@C=1.5pc{i\ar@/^.5pc/[r]^{(p,q)} & i\ar@/^.5pc/[l]^{(r,s)}},
\xymatrix@R=1.5pc@C=1.5pc{i\ar[r]^{(u,v)} & i}
\right\}
& = & 
-\delta_{pv}\xymatrix@R=1.5pc@C=1.5pc{i\ar@/^.5pc/[r]^{(u,q)} & i\ar@/^.5pc/[l]^{(r,s)}}
+\delta_{us}\xymatrix@R=1.5pc@C=1.5pc{i\ar@/^.5pc/[r]^{(p,q)} & i\ar@/^.5pc/[l]^{(r,v)}}
+\delta_{vr}\xymatrix@R=1.5pc@C=1.5pc{i\ar@/^.5pc/[r]^{(p,q)} & i\ar@/^.5pc/[l]^{(u,s)}}
-\delta_{uq}\xymatrix@R=1.5pc@C=1.5pc{i\ar@/^.5pc/[r]^{(p,v)} & i\ar@/^.5pc/[l]^{(r,s)}}
\end{eqnarray*}
\end{proof}
It becomes clear from the proof that the double Poisson-Lychnerowicz cohomology for a generic finite dimensional semi-simple algebra will never become zero. For the simplest non-trivial semi-simple algebra $S=\CC\oplus\CC$ one easily sees that
$$\dim H^i_P(S)  = \left\{\begin{array}{ll} 1 & i\neq 0 \mathrm{~even} \\ 0 & \mathrm{otherwise}\end{array}\right.$$

\section{Double Derivations on Amalgamated Products}\label{DDAP}
In \cite{JanLievenTSA}, the notion of a tree of semisimple algebras was introduced and linked to the representation theory of $\SL_2(\ZZ)$ via the amalgamated product $\ZZ_6*_{\ZZ_2}\ZZ_4$. The study of the finite dimensional representations of this amalgamated product can be seen as equivalent to the study of the finite dimensional representations of the amalgamated product of semisimple algebras $\CC^{\oplus 6}*_{\CC^{\oplus 2}}\CC^{\oplus 4}$. More generally, the representation theory of any torus knot group can be reduced via the representation theory of $\ZZ_p*\ZZ_q$ to the representation theory of $\CC^{\oplus p}*\CC^{\oplus q}$. Adriaenssens and Le~Bruyn show in \cite{JanLievenTSA} that the representation spaces of these amalgamated products have an \'etale cover characterized by a symmetric quiver and hence have a double Poisson structure as well. 

From \cite[Prop 2.4.1]{MichelDPA} we know that each pair of double Poisson brackets on the factors of an amalgamated product induces a unique double Poisson bracket on the product, so combining the nonzero double Poisson brackets we determined in the previous section of the paper, we obtain double Poisson structures on the representation spaces of the amalgamated products themselves rather than on an \'etale cover. In this section, we will formulate an explicit description of the Poisson brackets they induce on the quotient space $\iss_n(S*T)$

We begin by fixing notation for the rest of this section.
We let $S = \CC^{\oplus n}$ with orthogonal idempotents $e_1,\dots, e_n$ and $T = \CC^{\oplus m}$ with orthogonal idempotents $f_1, \dots, f_m$. By $Q_S$ and $Q_T$ we denote as before the  double derivation quivers of $S$ and $T$. The arrows in $Q_S$ will be denoted by $a_{ij}$ and the arrows in $Q_T$ will be $b_{kl}$ with $h(a_{ij}) = i$, $t(a_{ij}) = j$, $h(b_{kl}) = k$ and $t(b_{kl}) = l$. Finally, we let $P = \sum_{i,j} c_{ij}a_{ij}a_{ji}$ be a double Poisson tensor in $\CC Q_S$ and $P' = \sum_{k,l} d_{kl}b_{kl}b_{lk}$ a double Poissson tensor in $\CC Q_T$. We have the following expression for the double Poisson bracket induced by $P$ and $P'$ on $S*T$.
\begin{lemma}\label{DoubeStructureAmalgamated}
Let $x = e_{i_1}*f_{i_1'}*\dots *e_{i_p}*f_{i_p'}$ and $y = e_{j_1}*f_{j_1'}*\dots *e_{j_q}*f_{j_q'}$. Define for all $1\leq i,j \leq n$ and $1\leq k,\ell\leq m$ the elements $\overline{c}_{ij} = c_{ij} - c_{ji}$ and $\overline{d}_{k\ell} = d_{k\ell}-d_{\ell k}$ in $\CC$ and the elements $\overline{e}_i = -\sum_{s\neq i} \overline{c}_{is}e_s$ and $\overline{f}_k = -\sum_{s\neq i} \overline{d}_{ks}e_s$ then
\begin{align*}
\db{x,y} =
& \sum_{\ell = 1}^p 
   & \sum_{\stackrel{k=1}{j_k'\neq i_\ell'}}^q
      & \overline{d}_{i_\ell' j_k'}( e_{j_1}*\dots * e_{j_k}* (f_{i_\ell'}*\dots *f_{i_p'}\otimes e_{i_1}*\dots * e_{i_\ell}*f_{j_k'} \\
& & & 
- f_{j_k'}* e_{i_{\ell+1}}*\dots *f_{i_p'}\otimes e_{i_1}*\dots * f_{i_\ell'})*e_{j_{k+1}}*\dots * f_{j_q'}) \\
& & + \sum_{\stackrel{k=1}{j_k\neq i_\ell}}^q &
 \overline{c}_{i_\ell j_k} (e_{j_1}*\dots * f_{j_{k-1}'} *(e_{i_\ell}*\dots * f_{i_p'}\otimes e_{i_1}*\dots * f_{i_{\ell-1}'}* e_{j_k}\\
& &  &- e_{j_k}*f_{i_\ell'}*\dots * f_{i_p'}\otimes e_{i_1}*\dots * e_{i_\ell})*f_{j_k'}*\dots *f_{j_q'}) \\
& & + \sum_{\stackrel{k=1}{j_k' = i_\ell'}}^q
& e_{j_1}*\dots * e_{j_k}* (f_{i_\ell'}*\dots *f_{i_p'}\otimes e_{i_1}*\dots * e_{i_\ell}*\overline{f}_{j_k'} \\
& & &
- \overline{f}_{j_k'}* e_{i_{\ell+1}}*\dots *f_{i_p'}\otimes e_{i_1}*\dots * f_{i_\ell'})*e_{j_{k+1}}*\dots * f_{j_q'} \\
& &+ \sum_{\stackrel{k=1}{j_k = i_\ell}}^q &
e_{j_1}*\dots * f_{j_{k-1}'} *(e_{i_\ell}*\dots * f_{i_p'}\otimes e_{i_1}*\dots * f_{i_{\ell-1}'}* \overline{e}_{j_k}\\
& &  & - \overline{e}_{j_k}*f_{i_\ell'}*\dots * f_{i_p'}\otimes e_{i_1}*\dots * e_{i_\ell})*f_{j_k'}*\dots *f_{j_q'}.
\end{align*}
\end{lemma}
\begin{proof}
We will prove the statement by induction on the length $2p$ of $x$. First of all note that
$$\db{f_b, e_{j_1}*f_{j_1}*\dots *e_{j_q}*f_{j_q}} = \sum_{k=1}^q e_{j_1}*\dots * e_{j_k}\db{f_b,f_{j_k}}*e_{j_{k+1}}*\dots *f_{j_q}$$
and
$$\db{e_a, e_{j_1}*f_{j_1}*\dots *e_{j_q}*f_{j_q}} = \sum_{k=1}^{q} e_{j_1}*\dots * f_{j_{k-1}}\db{e_a,e_{j_{k}}}*f_{j_{k}}*\dots *f_{j_q}$$
because $\db{f_b,-}$ is a $S$-linear double derivation and $\db{e_a,-}$ is a $T$-linear double derivation. This means for $p=1$ that
\begin{eqnarray*}
\db{e_a*f_b,y} &=& e_a\circ\db{f_b,y} + \db{e_a,y}\circ f_b \\
& = & e_a\circ\sum_{k=1}^q e_{j_1}*\dots * e_{j_k}*\db{f_b,f_{j_k}}*e_{j_{k+1}}*\dots *f_{j_q} \\
&    & + \sum_{k=1}^{q} e_{j_1}*\dots * f_{j_{k-1}}\db{e_a,e_{j_{k}}}*f_{j_{k}}*\dots *f_{j_q})\circ f_b.
\end{eqnarray*}
Now by definition of the double brackets associated to $P$ and $P'$ we get
\begin{align*}
\db{f_b,f_{j_k}} & = (d_{bj_k} - d_{j_kb})(f_b\otimes f_{j_k} - f_{j_k}\otimes f_b) & b\neq j_k \\
\db{f_b,f_b}& = -\sum_{s\neq b} (d_{bs} - d_{sb})(f_b\otimes f_s - f_s\otimes f_b) \\
& = f_b\otimes \overline{f}_b - \overline{f}_b\otimes f_b
\end{align*}
and
\begin{align*}
\db{e_a,e_{j_k}}& =  (c_{aj_k} - c_{j_ka})(e_a\otimes e_{j_k} - e_{j_k}\otimes e_a) & a\neq j_k\\
\db{e_a,e_a}& = -\sum_{s\neq a} (c_{as} - c_{sa})(e_a\otimes e_{s} - e_{s}\otimes e_a)\\
&=e_a\otimes\overline{e}_a - \overline{e}_a\otimes e_a.
\end{align*}
So
\begin{eqnarray*}
\db{e_a*f_b,y} & = & \sum_{\stackrel{k=1}{j_k'\neq b}}^q \overline{d}_{bj_k}e_{j_1}*\dots * e_{j_k}*(f_b\otimes e_a*f_{j_k} - f_{j_k}\otimes e_a*f_b)*e_{j_{k+1}}*\dots *f_{j_q} \\
&  &+ \sum_{\stackrel{k=1}{j_k\neq a}}^{q} \overline{c}_{aj_k}e_{j_1}*\dots * f_{j_{k-1}'}*(e_a*f_b\otimes e_{j_k} - e_{j_k}*f_b\otimes e_a)*f_{j_{k}'}*\dots *f_{j_q'} \\
& & + \sum_{\stackrel{k=1}{j_k' = b}}^q e_{j_1}*\dots * e_{j_k}*(f_b\otimes e_a*\overline{f}_{b} - \overline{f}_{b}\otimes e_a*f_b)*e_{j_{k+1}}*\dots *f_{j_q'} \\
&  &+ \sum_{\stackrel{k=1}{j_k = a}}^{q} e_{j_1}*\dots * f_{j_{k-1}'}*(e_a*f_b\otimes \overline{e}_a -\overline{e}_a*f_b\otimes e_a)*f_{j_{k}'}*\dots *f_{j_q'} 
\end{eqnarray*}
which means the claim holds for $p=1$. Now assume the lemma is correct up to $p-1$, then
$$\db{x,y} = (e_{i_1}*\dots * f_{i_{p-1}'})\circ\db{e_{i_p}*f_{i_p'},y} + \db{e_{i_1}*\dots * f_{i_{p-1}'},y}\circ (e_{i_p}*f_{i_p'}).
$$
The first term in this expression is equal to
\begin{align*}
\sum_{\stackrel{k=1}{j_k'\neq i_\ell'}}^q
& \overline{d}_{i_p j_k} (e_{j_1}*\dots * e_{j_k}* (f_{i_p'}*\dots *f_{i_p'}\otimes e_{i_1}*\dots * e_{i_p}*f_{j_k'}
\\
& 
- f_{j_k'}* e_{i_{p+1}}*\dots *f_{i_p'}\otimes e_{i_1}*\dots * f_{i_p'})*e_{j_{k+1}}*\dots * f_{j_q'})
 \\
+ \sum_{\stackrel{k=1}{j_k\neq i_p}}^q
&
\overline{c}_{i_p j_k}( e_{j_1}*\dots * f_{j_{k-1}'} *(e_{i_p}*\dots * f_{i_p'}\otimes e_{i_1}*\dots * f_{i_{p-1}'}* e_{j_k}
 \\
&- e_{j_k}*f_{i_p'}*\dots * f_{i_p'}\otimes e_{i_1}*\dots * e_{i_p})*f_{j_k'}*\dots *f_{j_q'})
\\
+ \sum_{\stackrel{k=1}{j_k' = i_p'}}^q
& e_{j_1}*\dots * e_{j_k}* (f_{i_p'}*\dots *f_{i_p'}\otimes e_{i_1}*\dots * e_{i_p}*\overline{f}_{j_k'}
\\
& 
- \overline{f}_{j_k'}* e_{i_{p+1}}*\dots *f_{i_p'}\otimes e_{i_1}*\dots * f_{i_p'})*e_{j_{k+1}}*\dots * f_{j_q'} \\
+ \sum_{\stackrel{k=1}{j_k = i_p}}^q &
e_{j_1}*\dots * f_{j_{k-1}'} *(e_{i_p}*\dots * f_{i_p'}\otimes e_{i_1}*\dots * f_{i_{p-1}'}* \overline{e}_{j_k}
\\
& - \overline{e}_{j_k}*f_{i_p'}*\dots * f_{i_p'}\otimes e_{i_1}*\dots * e_{i_p})*f_{j_k'}*\dots *f_{j_q'}.
\end{align*}
and the second term becomes
\begin{align*}
\db{x,y} =
& \sum_{\ell = 1}^{p-1}
   & \sum_{\stackrel{k=1}{j_k'\neq i_\ell'}}^q
      & \overline{d}_{i_\ell j_k} (e_{j_1}*\dots * e_{j_k}* (f_{i_\ell'}*\dots *f_{i_p'}\otimes e_{i_1}*\dots * e_{i_\ell}*f_{j_k'} \\
& & & 
- f_{j_k'}* e_{i_{\ell+1}}*\dots *f_{i_p'}\otimes e_{i_1}*\dots * f_{i_\ell'})*e_{j_{k+1}}*\dots * f_{j_q'}) \\
& & + \sum_{\stackrel{k=1}{j_k\neq i_\ell}}^q &
 \overline{c}_{i_\ell j_k}( e_{j_1}*\dots * f_{j_{k-1}'} *(e_{i_\ell}*\dots * f_{i_p'}\otimes e_{i_1}*\dots * f_{i_{\ell-1}'}* e_{j_k}\\
& &  &- e_{j_k}*f_{i_\ell'}*\dots * f_{i_p'}\otimes e_{i_1}*\dots * e_{i_\ell})*f_{j_k'}*\dots *f_{j_q'}) \\
& & + \sum_{\stackrel{k=1}{j_k' = i_\ell'}}^q
& e_{j_1}*\dots * e_{j_k}* (f_{i_\ell'}*\dots *f_{i_p'}\otimes e_{i_1}*\dots * e_{i_\ell}*\overline{f}_{j_k'} \\
& & &
- \overline{f}_{j_k'}* e_{i_{\ell+1}}*\dots *f_{i_p'}\otimes e_{i_1}*\dots * f_{i_\ell'})*e_{j_{k+1}}*\dots * f_{j_q'} \\
& &+ \sum_{\stackrel{k=1}{j_k = i_\ell}}^q &
e_{j_1}*\dots * f_{j_{k-1}'} *(e_{i_\ell}*\dots * f_{i_p'}\otimes e_{i_1}*\dots * f_{i_{\ell-1}'}* \overline{e}_{j_k}\\
& &  & - \overline{e}_{j_k}*f_{i_\ell'}*\dots * f_{i_p'}\otimes e_{i_1}*\dots * e_{i_\ell})*f_{j_k'}*\dots *f_{j_q'}.
\end{align*}
Adding these two terms together then yields the expression from the lemma.
\end{proof}

\begin{prop}\label{InducedStructureAmalgamated}
With notations as before, the Poisson bracket induced by $P$ and $P'$ on $\iss(S*T)$ is equal to
\begin{eqnarray*}
\{tr(x),tr(y)\} &=& 
\sum_{\ell = 1}^p 
tr({\sigma^{2\ell-1}(x)}*
(\sum_{\stackrel{k=1}{j_k'\neq i_\ell'}}^q \overline{d}_{i_\ell' j_k'}{\sigma^{2k-1}(y)} - \sum_{\stackrel{k=1}{j_k\neq i_\ell}}^q  \overline{c}_{i_\ell j_k}{\sigma^{2k-1}(y)})) \\
& & 
-tr({\sigma^{2\ell}(x)}*(\sum_{\stackrel{k=1}{j_k'\neq i_\ell'}}^q \overline{d}_{i_\ell' j_k'}{\sigma^{2k}(y)} - \sum_{\stackrel{k=1}{j_k\neq i_\ell}}^q  \overline{c}_{i_{\ell+1} j_k}{\sigma^{2k}(y)}) ) \\
& & + \sum_{\stackrel{k=1}{j_k' = i_\ell'}}^q
tr({\sigma^{2\ell-1}(x)}*{\sigma^{2k-1}(y_k)})
- tr({\sigma^{2\ell}(x)}*{\sigma^{2k}(y_k)}) \\ 
& & 
+ \sum_{\stackrel{k=1}{j_k = i_\ell}}^q
tr({\sigma^{2(\ell-1)}(x)}*{\sigma^{2(k-1)}({}_ky)})
- tr({\sigma^{2\ell-1}(x)}*{\sigma^{2k-1}({}_ky)}).
\end{eqnarray*}
with $\sigma(a*x) = x*a$ the cyclic permutation and where $y_k = e_{j_1}*f_{j_1'}*\dots *f_{j_{k-1}'}* e_{j_k}*\overline{f}_{j_k}*e_{j_{k+1}}*\dots *e_{j_q}*f_{j_q'}$ and ${}_ky = e_{j_1}*f_{j_1'}*\dots *f_{j_{k-1}'}* \overline{e}_{j_k}*{f}_{j_k}*\dots *e_{j_q}*f_{j_q'}$.
\end{prop}
\begin{proof}
By definition of the bracket induced on $\iss_n(S*T)$ we have
$$\{ tr(x), tr(y)\} = \db{x,y}_{ij}\db{x,y}_{ji}.$$
Using the expression found in the previous lemma this becomes
\begin{align*}
\sum_{\ell = 1}^p 
   & \sum_{\stackrel{k=1}{j_k'\neq i_\ell'}}^q
      & 
\overline{d}_{i_\ell j_k}(
(e_{j_1}*\dots * e_{j_k}*f_{i_\ell'}*\dots *f_{i_p'})_{rs}(e_{i_1}*\dots * e_{i_\ell}*f_{j_k'}*e_{j_{k+1}}*\dots * f_{j_q'})_{sr} \\
& & 
- (e_{j_1}*\dots * e_{j_k}*f_{j_k'}* e_{i_{\ell+1}}*\dots *f_{i_p'})_{rs}(e_{i_1}*\dots * f_{i_\ell'}*e_{j_{k+1}}*\dots * f_{j_q'})_{sr} )\\
& + \sum_{\stackrel{k=1}{j_k\neq i_\ell}}^q &
 \overline{c}_{i_\ell j_k}(
 (e_{j_1}*\dots * f_{j_{k-1}'} *e_{i_\ell}*\dots * f_{i_p'})_{rs}(e_{i_1}*\dots * f_{i_{\ell-1}'}* e_{j_k}*f_{j_k'}*\dots *f_{j_q'})_{sr}\\
&  &- (e_{j_1}*\dots * f_{j_{k-1}'} *e_{j_k}*f_{i_\ell'}*\dots * f_{i_p'})_{rs}(e_{i_1}*\dots * e_{i_\ell}*f_{j_k'}*\dots *f_{j_q'})_{sr} )\\
& + \sum_{\stackrel{k=1}{j_k' = i_\ell'}}^q
& 
(e_{j_1}*\dots * e_{j_k}*f_{i_\ell'}*\dots *f_{i_p'})_{rs}(e_{i_1}*\dots * e_{i_\ell}*\overline{f}_{j_k'}*e_{j_{k+1}}*\dots * f_{j_q'})_{sr} \\
& &
- (e_{j_1}*\dots * e_{j_k}*\overline{f}_{j_k'}* e_{i_{\ell+1}}*\dots *f_{i_p'})_{rs}(e_{i_1}*\dots * f_{i_\ell'}*e_{j_{k+1}}*\dots * f_{j_q'})_{sr}\\
&+ \sum_{\stackrel{k=1}{j_k = i_\ell}}^q &
(e_{j_1}*\dots * f_{j_{k-1}'} *e_{i_\ell}*\dots * f_{i_p'})_{rs}(e_{i_1}*\dots * f_{i_{\ell-1}'}* \overline{e}_{j_k}*f_{j_k'}*\dots *f_{j_q'})_{sr}\\
&  & - (e_{j_1}*\dots * f_{j_{k-1}'} *\overline{e}_{j_k}*f_{i_\ell'}*\dots * f_{i_p'})_{rs}(e_{i_1}*\dots * e_{i_\ell}*f_{j_k'}*\dots *f_{j_q'})_{sr}.
\end{align*}
This is equal to
\begin{align*}
\sum_{\ell = 1}^p 
   & \sum_{\stackrel{k=1}{j_k'\neq i_\ell'}}^q
      & 
\overline{d}_{i_\ell j_k}(
tr(e_{j_1}*\dots * e_{j_k}*f_{i_\ell'}*\dots *f_{i_p'}*e_{i_1}*\dots * e_{i_\ell}*f_{j_k'}*e_{j_{k+1}}*\dots * f_{j_q'}) \\
& & 
- tr(e_{j_1}*\dots * e_{j_k}*f_{j_k'}* e_{i_{\ell+1}}*\dots *f_{i_p'}*e_{i_1}*\dots * f_{i_\ell'}*e_{j_{k+1}}*\dots * f_{j_q'}))\\
& + \sum_{\stackrel{k=1}{j_k\neq i_\ell}}^q &
 \overline{c}_{i_\ell j_k}(
 tr(e_{j_1}*\dots * f_{j_{k-1}'} *e_{i_\ell}*\dots * f_{i_p'}*e_{i_1}*\dots * f_{i_{\ell-1}'}* e_{j_k}*f_{j_k'}*\dots *f_{j_q'})\\
&  &- tr(e_{j_1}*\dots * f_{j_{k-1}'} *e_{j_k}*f_{i_\ell'}*\dots * f_{i_p'}*e_{i_1}*\dots * e_{i_\ell}*f_{j_k'}*\dots *f_{j_q'}))\\
& + \sum_{\stackrel{k=1}{j_k' = i_\ell'}}^q
& 
tr(e_{j_1}*\dots * e_{j_k}*f_{i_\ell'}*\dots *f_{i_p'}*e_{i_1}*\dots * e_{i_\ell}*\overline{f}_{j_k'}*e_{j_{k+1}}*\dots * f_{j_q'}) \\
& &
- tr(e_{j_1}*\dots * e_{j_k}*\overline{f}_{j_k'}* e_{i_{\ell+1}}*\dots *f_{i_p'}*e_{i_1}*\dots * f_{i_\ell'}*e_{j_{k+1}}*\dots * f_{j_q'})\\
&+ \sum_{\stackrel{k=1}{j_k = i_\ell}}^q &
tr(e_{j_1}*\dots * f_{j_{k-1}'} *e_{i_\ell}*\dots * f_{i_p'}*e_{i_1}*\dots * f_{i_{\ell-1}'}* \overline{e}_{j_k}*f_{j_k'}*\dots *f_{j_q'})\\
&  & - tr(e_{j_1}*\dots * f_{j_{k-1}'} *\overline{e}_{j_k}*f_{i_\ell'}*\dots * f_{i_p'}*e_{i_1}*\dots * e_{i_\ell}*f_{j_k'}*\dots *f_{j_q'}).
\end{align*}
Using the fact that traces are invariant under cyclic permutation, we may rewrite this as
\begin{align*}
\sum_{\ell = 1}^p 
   & \sum_{\stackrel{k=1}{j_k'\neq i_\ell'}}^q
      & 
\overline{d}_{i_\ell j_k}(
tr(
\underbrace{f_{i_\ell'}*\dots *f_{i_p'}*e_{i_1}*\dots * e_{i_\ell}}_{\sigma^{2\ell-1}(x)}
*
\underbrace{f_{j_k'}*e_{j_{k+1}}*\dots * f_{j_q'}*e_{j_1}*\dots * e_{j_k})}_{\sigma^{2k-1}(y)} \\
& & 
- tr(
\underbrace{e_{i_{\ell+1}}*\dots *f_{i_p'}*e_{i_1}*\dots * f_{i_\ell'}}_{\sigma^{2\ell}(x)}
*
\underbrace{e_{j_{k+1}}*\dots * f_{j_q'}*e_{j_1}*\dots * e_{j_k}*f_{j_k'}}_{\sigma^{2k}(y)}))\\
& + \sum_{\stackrel{k=1}{j_k\neq i_\ell}}^q &
 \overline{c}_{i_\ell j_k}(
tr(
\underbrace{e_{i_\ell}*\dots * f_{i_p'}*e_{i_1}*\dots * f_{i_{\ell-1}'}}_{\sigma^{2(\ell - 1)}(x)}
*
\underbrace{e_{j_k}*f_{j_k'}*\dots *f_{j_q'}*e_{j_1}*\dots * f_{j_{k-1}'}}_{\sigma^{2(k-1)}(y)})\\
&  &
- tr(
\underbrace{*f_{i_\ell'}*\dots * f_{i_p'}*e_{i_1}*\dots * e_{i_\ell}}_{\sigma^{2\ell-1}(x)}
*
\underbrace{f_{j_k'}*\dots *f_{j_q'}*e_{j_1}*\dots * f_{j_{k-1}'}*e_{j_k}}_{\sigma^{2k-1}(y)}))\\
& + \sum_{\stackrel{k=1}{j_k' = i_\ell'}}^q
& 
tr(
\underbrace{f_{i_\ell'}*\dots *f_{i_p'}*e_{i_1}*\dots * e_{i_\ell}}_{\sigma^{2\ell-1}(x)}
*
\underbrace{\overline{f}_{j_k'}*e_{j_{k+1}}*\dots * f_{j_q'}*e_{j_1}*\dots * e_{j_k}}_{\sigma^{2k-1}(y_k)}) \\
& &
- tr(
\underbrace{e_{i_{\ell+1}}*\dots *f_{i_p'}*e_{i_1}*\dots * f_{i_\ell'}}_{\sigma^{2\ell}(x)}
*
\underbrace{e_{j_{k+1}}*\dots * f_{j_q'}*e_{j_1}*\dots * e_{j_k}*\overline{f}_{j_k'}}_{\sigma^{2k}(y_k)})\\
&+ \sum_{\stackrel{k=1}{j_k = i_\ell}}^q &
tr(
\underbrace{e_{i_\ell}*\dots * f_{i_p'}*e_{i_1}*\dots * f_{i_{\ell-1}'}}_{\sigma^{2(\ell-1)}(x)}
*
\underbrace{\overline{e}_{j_k}*f_{j_k'}*\dots *f_{j_q'}*e_{j_1}*\dots * f_{j_{k-1}'}}_{\sigma^{2(k-1)}({}_ky)})\\
&  & 
- tr(
\underbrace{f_{i_\ell'}*\dots * f_{i_p'}*e_{i_1}*\dots * e_{i_\ell}}_{\sigma^{2\ell-1}(x)}
*
\underbrace{f_{j_k'}*\dots *f_{j_q'}*e_{j_1}*\dots * f_{j_{k-1}'} *\overline{e}_{j_k}}_{\sigma^{2k-1}({}_ky)}).
\end{align*}
Regrouping this expression then gives the expression from the proposition.
\end{proof}
Although this expression yields nontrivial Poisson brackets in general, note that we have the following corollary.
\begin{cor}
For $S = \CC^{\oplus 2}$ and $T = \CC^{\oplus d}$ with $d$ arbitrary, any Poisson bracket induced by double Poisson tensors $P$ and $P'$ becomes the zero-bracket.
\end{cor}
\begin{proof}
First of all note that for $S = \CC^{\oplus 2}$ the ring of invariants $\CC[\rep_n(S*T)]^{\GL_n}$ is generated by expressions of the form $tr(e_1*f_{i_1}*e_1*f_{i_2}\dots *e_1*f_{i_p})$. But then
\begin{eqnarray*}
\sigma^{2\ell-1}(x)*\sigma^{2k-1}(y) &=& \sigma^{2\ell}(x)\sigma^{2k}(y)\\
{\sigma^{2\ell-1}(x)}*{\sigma^{2k-1}(y_k)} & = & {\sigma^{2\ell}(x)}*{\sigma^{2k}(y_k)} \\ 
{\sigma^{2(\ell-1)}(x)}*{\sigma^{2(k-1)}({}_ky)} & = & {\sigma^{2\ell-1}(x)}*{\sigma^{2k-1}({}_ky)}
\end{eqnarray*}
up to cyclic permutation, so the bracket becomes zero.
\end{proof}

From Theorem \ref{DoubleDerivationsRelative} we know the module of $\CC^{\oplus 2}$-relative double derivations on $\CC^{\oplus 4}$ can be depicted as
$$\xymatrix@R=1.5em@C=1.5em{
\vtx{~}\ar@/^.5em/[r] & \vtx{~}\ar@/^.5em/[l] & \vtx{~}\ar@/^.5em/[r] & \vtx{~}\ar@/^.5em/[l]
}$$
and the module of $\CC^{\oplus 2}$-relative double derivations on $\CC^{\oplus 6}$ can be depicted as 
$$
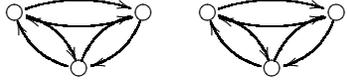
\xymatrix@R=1.5em@C=1.5em{
\vtx{~}\ar@/^.5em/[rr]\ar@/^.5em/[dr] & & \vtx{~}\ar@/^.5em/[ll]\ar@/^.5em/[dl] & \vtx{~}\ar@/^.5em/[rr]\ar@/^.5em/[dr] &  & \vtx{~}\ar@/^.5em/[ll]\ar@/^.5em/[dl] \\
& \vtx{~}\ar@/^.5em/[ur]\ar@/^.5em/[ul] & &  & \vtx{~}\ar@/^.5em/[ur]\ar@/^.5em/[ul]
}.$$
The last corollary then also means that we only get the trivial bracket on this amalgamated product, so in case of $\SL_2(\ZZ)$ we do not obtain a Poisson structure on the quotient spaces from the double Poisson structures on the factors, whereas Adriaenssens and Le~Bruyn obtain nontrivial Poisson structures on an \'etale cover. For higher dimensional factors, however, the above formula yields nontrivial Poisson brackets on the quotient variety.

\end{document}